\documentclass[10pt]{amsart}
\usepackage{latexsym}
\usepackage[dvips]{color}
\usepackage{float}
\usepackage{afterpage} 
\usepackage [pdftex]{graphicx}
\usepackage{epsfig}
\usepackage{amsmath}
\usepackage{amsthm}
\usepackage{amsfonts}
\usepackage{amssymb}
\usepackage{setspace}
\usepackage{subfigure}
\usepackage{tikz}
\usetikzlibrary{arrows}
\usetikzlibrary{arrows.meta}
\usetikzlibrary{matrix}
\usepackage{subfigure}

\usepackage{colortbl}
\usepackage{url}

\newtheorem{definition}{Definition}
\newtheorem{proposition}{Proposition}

\newtheorem{theorem}{Theorem}

\begin{document}

\newcommand{\autdiagrampositive}[8] {
\scalebox{0.6}{
\begin{tikzpicture}[scale=1.0]
\node (Aug0_top) at (0,0) {$C_{\text{aug}}$};
\node (Aug3_top) at (1,0) {$G_{\text{aug}}$};
\node (Aug6_top) at (2,0) {$D_{\text{aug}}$};
\node (Aug9_top) at (3,0) {$F_{\text{aug}}$};

\node (Aug0_bottom) at (0,-2.5) {$C_{\text{aug}}$};
\node (Aug3_bottom) at (1,-2.5) {$G_{\text{aug}}$};
\node (Aug6_bottom) at (2,-2.5) {$D_{\text{aug}}$};
\node (Aug9_bottom) at (3,-2.5) {$F_{\text{aug}}$};

\draw[-,>=latex, line width=1.0] (Aug0_top) to [in=90,out=270] (Aug#1_bottom) ;
\draw[-,>=latex, line width=1.0] (Aug3_top) to [in=90,out=270] (Aug#2_bottom) ;
\draw[-,>=latex, line width=1.0] (Aug6_top) to [in=90,out=270] (Aug#3_bottom) ;
\draw[-,>=latex, line width=1.0] (Aug9_top) to [in=90,out=270] (Aug#4_bottom) ;

\node (Top_0) at (5,0) {$\overline{C_M}$};
\node (Top_3) at (6,0) {$\overline{G_M}$};
\node (Top_6) at (7,0) {$\overline{D_M}$};
\node (Top_9) at (8,0) {$\overline{F_M}$};

\node (Bottom_0) at (5,-2.5) {$\overline{C_M}$};
\node (Bottom_3) at (6,-2.5) {$\overline{G_M}$};
\node (Bottom_6) at (7,-2.5) {$\overline{D_M}$};
\node (Bottom_9) at (8,-2.5) {$\overline{F_M}$};

\draw[-,>=latex, line width=1.0] (Top_0) to [in=90,out=270] node [midway,fill=white,text=blue,pos=0.85]  {$g_0$} (Bottom_#5)  ;
\draw[-,>=latex, line width=1.0] (Top_3) to [in=90,out=270] node [midway,fill=white,text=blue,pos=0.85]  {$g_1$} (Bottom_#6) ;
\draw[-,>=latex, line width=1.0] (Top_6) to [in=90,out=270] node [midway,fill=white,text=blue,pos=0.85]  {$g_2$} (Bottom_#7) ;
\draw[-,>=latex, line width=1.0] (Top_9) to [in=90,out=270] node [midway,fill=white,text=blue,pos=0.85]  {$g_3$} (Bottom_#8) ;
\end{tikzpicture}
}
}

\newcommand{\autdiagramnegative}[8] {
\scalebox{0.6}{
\begin{tikzpicture}[scale=1.0]
\node (Aug0_top) at (0,0) {$C_{\text{aug}}$};
\node (Aug3_top) at (1,0) {$G_{\text{aug}}$};
\node (Aug6_top) at (2,0) {$D_{\text{aug}}$};
\node (Aug9_top) at (3,0) {$F_{\text{aug}}$};

\node (Aug0_bottom) at (0,-2.5) {$C_{\text{aug}}$};
\node (Aug3_bottom) at (1,-2.5) {$G_{\text{aug}}$};
\node (Aug6_bottom) at (2,-2.5) {$D_{\text{aug}}$};
\node (Aug9_bottom) at (3,-2.5) {$F_{\text{aug}}$};

\draw[-,>=latex, line width=1.0] (Aug0_top) to [in=90,out=270] (Aug#1_bottom) ;
\draw[-,>=latex, line width=1.0] (Aug3_top) to [in=90,out=270] (Aug#2_bottom) ;
\draw[-,>=latex, line width=1.0] (Aug6_top) to [in=90,out=270] (Aug#3_bottom) ;
\draw[-,>=latex, line width=1.0] (Aug9_top) to [in=90,out=270] (Aug#4_bottom) ;

\node (Top_0) at (5,0) {$\overline{C_M}$};
\node (Top_3) at (6,0) {$\overline{G_M}$};
\node (Top_6) at (7,0) {$\overline{D_M}$};
\node (Top_9) at (8,0) {$\overline{F_M}$};

\node (Bottom_0) at (5,-2.5) {$\overline{C_m}$};
\node (Bottom_3) at (6,-2.5) {$\overline{G_m}$};
\node (Bottom_6) at (7,-2.5) {$\overline{D_m}$};
\node (Bottom_9) at (8,-2.5) {$\overline{F_m}$};

\draw[-,>=latex, line width=1.0] (Top_0) to [in=90,out=270] node [midway,fill=white,text=blue,pos=0.85]  {$g_0$} (Bottom_#5)  ;
\draw[-,>=latex, line width=1.0] (Top_3) to [in=90,out=270] node [midway,fill=white,text=blue,pos=0.85]  {$g_1$} (Bottom_#6) ;
\draw[-,>=latex, line width=1.0] (Top_6) to [in=90,out=270] node [midway,fill=white,text=blue,pos=0.85]  {$g_2$} (Bottom_#7) ;
\draw[-,>=latex, line width=1.0] (Top_9) to [in=90,out=270] node [midway,fill=white,text=blue,pos=0.85]  {$g_3$} (Bottom_#8) ;
\end{tikzpicture}
}
}

\title[Composing (with) automorphisms in the colored Cube Dance]{Composing (with) automorphisms in the colored Cube Dance: an interactive tool for musical chord transformation.}
\author{Alexandre Popoff}
\address{Independent Researcher}
\email{al.popoff@free.fr}

\author{Corentin Guichaoua}
\address{STMS (UMR9912), CNRS, Ircam, Sorbonne Universit\'e, Minist\`ere de la Culture, Paris, France}
\email{Corentin.Guichaoua@ircam.fr}

\author{Moreno Andreatta}
\address{CNRS/Institute for Advanced Mathematical Research, ITI CREAA, University of Strasbourg, France and IRCAM, Paris, France}
\email{andreatta@math.unistra.fr, andreatta@ircam.fr}


\begin{abstract}
The `colored Cube Dance' is an extension of Douthett's and Steinbach's Cube Dance graph, related to a monoid of binary relations defined on the set of major, minor, and augmented triads. This contribution explores the automorphism group of this monoid action, as a way to transform chord progressions. We show that this automorphism group is of order 7776 and is isomorphic to $({\mathbb{Z}_3}^4 \rtimes D_8) \rtimes (D_6 \times \mathbb{Z}_2)$. The size and complexity of this group makes it unwieldy: we therefore provide an interactive tool \textit{via} a web interface based on common HTML/Javascript frameworks for students, musicians, and composers to explore these automorphisms, showing the potential of these technologies for math/music outreach activities. 
\end{abstract}

\maketitle
\section{An algebraic introduction to the colored Cube Dance}

The Cube Dance is a well-known structure introduced by Douthett and Steinbach in their work on parsimonious graphs between triads \cite{Douthett1998}. Its definition involves the $\mathcal{P}_{1,0}$ binary relation which relates two pitch-class sets if they differ by a single pitch class a semitone apart. The Cube Dance is then defined as the graph having the major, minor, and augmented triads as its vertices and the set of pairs of triads related by $\mathcal{P}_{1,0}$ as its set of edges. Since the the classical neo-Riemannian $P$ and $L$ operations imply the $\mathcal{P}_{1,0}$ binary relation, some recent work \cite{Popoff2018} has investigated an extension of the Cube Dance wherein further refinements of the $\mathcal{P}_{1,0}$ relation are considered. More precisely, three binary relations $\mathcal{U}$, $\mathcal{P}$, and $\mathcal{L}$ are defined on the set of major, minor, and augmented triads as follows. The notation we adopt for these triads is of the form $x_{\text{s}}$, where $x$ is a pitch class (the root for major, and minor triads, or any note for augmented triads), and $s$ is a subscript (`M', `m', or `aug') indicating the type of triad.

\begin{definition}
Let $X$ be the set of the 24 major and minor triads and the four augmented triads.
\begin{itemize}
\item{The relation $\mathcal{P}$ is the symmetric relation which coincides with the neo-Riemannian $P$ operation on major and minor triads and is the identity relation on augmented triads.}
\item{The relation $\mathcal{L}$ is the symmetric relation which coincides with the neo-Riemannian $L$ operation on major and minor triads and is the identity relation on augmented triads.}
\item{The relation $\mathcal{U}$ is the symmetric relation which relates an augmented triad with a major or minor triad if they are related by the $\mathcal{P}_{1,0}$ relation.}
\end{itemize}
\end{definition}

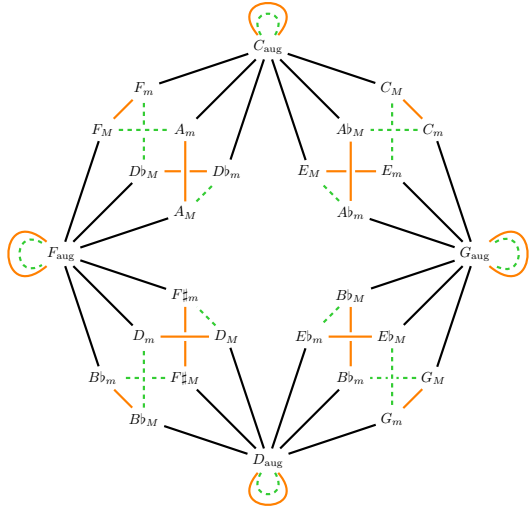
\begin{figure}[t!]
\begin{center}
\scalebox{0.55}{
\begin{tikzpicture}
\definecolor{color_U}{rgb}{0.0,0.0,0.0}
\definecolor{color_P}{rgb}{1.0,0.5,0.0}
\definecolor{color_L}{rgb}{0.2,0.8,0.2}

	\node (CM) at (3,4) {$C_M$};
	\node (Cm) at (4,3) {$C_m$};
	\node (Em) at (3,2) {$E_m$};
	\node (EM) at (1,2) {$E_M$};
	\node (Abm) at (2,1) {$A\flat_m$};
	\node (AbM) at (2,3) {$A\flat_M$};
	
	\draw[-,>=latex, line width=1.5, color_P] (EM) to (Em) ;
	\draw[-,>=latex, line width=1.5, color_L, dashed] (Em) to (CM) ;
	\draw[-,>=latex, line width=1.5, color_P] (CM) to (Cm) ;
		\draw[-,>=latex, line width=6, color=white] (Cm) to (AbM) ;
	\draw[-,>=latex, line width=1.5, color_L, dashed] (Cm) to (AbM) ;
		\draw[-,>=latex, line width=6, white] (AbM) to (Abm) ;
	\draw[-,>=latex, line width=1.5, color_P] (AbM) to (Abm) ;
	\draw[-,>=latex, line width=1.5, color_L, dashed] (Abm) to (EM) ;

	\node (Fm) at (-3,4) {$F_m$};
	\node (FM) at (-4,3) {$F_M$};
	\node (DbM) at (-3,2) {$D\flat_M$};
	\node (Dbm) at (-1,2) {$D\flat_m$};
	\node (AM) at (-2,1) {$A_M$};
	\node (Am) at (-2,3) {$A_m$};
	\draw[-,>=latex, line width=1.5, color_P] (Dbm) to (DbM) ;
	\draw[-,>=latex, line width=1.5, color_L, dashed] (DbM) to (Fm) ;
	\draw[-,>=latex, line width=1.5, color_P] (Fm) to (FM) ;
		\draw[-,>=latex, line width=6, white] (FM) to (Am) ;
	\draw[-,>=latex, line width=1.5, color_L, dashed] (FM) to (Am) ;
		\draw[-,>=latex, line width=6, white] (Am) to (AM) ;
	\draw[-,>=latex, line width=1.5, color_P] (Am) to (AM) ;
	\draw[-,>=latex, line width=1.5, color_L, dashed] (AM) to (Dbm) ;

	\node (BbM) at (-3,-4) {$B\flat_M$};
	\node (Bbm) at (-4,-3) {$B\flat_m$};
	\node (Dm) at (-3,-2) {$D_m$};
	\node (DM) at (-1,-2) {$D_M$};
	\node (Fsm) at (-2,-1) {$F\sharp_m$};
	\node (FsM) at (-2,-3) {$F\sharp_M$};
	\draw[-,>=latex, line width=1.5, color_P] (Fsm) to (FsM) ;
	\draw[-,>=latex, line width=1.5, color_L, dashed] (FsM) to (Bbm) ;
	\draw[-,>=latex, line width=1.5, color_P] (Bbm) to (BbM) ;
		\draw[-,>=latex, line width=6, white] (BbM) to (Dm) ;
	\draw[-,>=latex, line width=1.5, color_L, dashed] (BbM) to (Dm) ;
		\draw[-,>=latex, line width=6, white] (Dm) to (DM) ;
	\draw[-,>=latex, line width=1.5, color_P] (Dm) to (DM) ;
	\draw[-,>=latex, line width=1.5, color_L, dashed] (DM) to (Fsm) ;

	\node (Gm) at (3,-4) {$G_m$};
	\node (GM) at (4,-3) {$G_M$};
	\node (EbM) at (3,-2) {$E\flat_M$};
	\node (Ebm) at (1,-2) {$E\flat_m$};
	\node (BM) at (2,-1) {$B\flat_M$};
	\node (Bm) at (2,-3) {$B\flat_m$};
	\draw[-,>=latex, line width=1.5, color_L, dashed] (GM) to (Bm) ;
	\draw[-,>=latex, line width=1.5, color_P] (Bm) to (BM) ;
	\draw[-,>=latex, line width=1.5, color_L, dashed] (BM) to (Ebm) ;
		\draw[-,>=latex, line width=6, white] (Ebm) to (EbM) ;
	\draw[-,>=latex, line width=1.5, color_P] (Ebm) to (EbM) ;
		\draw[-,>=latex, line width=6, white] (EbM) to (Gm) ;
	\draw[-,>=latex, line width=1.5, color_L, dashed] (EbM) to (Gm) ;
	\draw[-,>=latex, line width=1.5, color_P] (Gm) to (GM) ;

	\node (Faug) at (-5,0) {$F_{\text{aug}}$};
	\node (Gaug) at (5,0) {$G_{\text{aug}}$};
	\node (Caug) at (0,5) {$C_{\text{aug}}$};
	\node (Daug) at (0,-5) {$D_{\text{aug}}$};

	\draw[-,>=latex, line width=1.5, color_U] (Fsm) to (Faug) ;
	\draw[-,>=latex, line width=1.5, color_U] (Dm) to (Faug) ;
	\draw[-,>=latex, line width=1.5, color_U] (Bbm) to (Faug) ;
	\draw[-,>=latex, line width=1.5, color_U] (FM) to (Faug) ;
	\draw[-,>=latex, line width=1.5, color_U] (DbM) to (Faug) ;
	\draw[-,>=latex, line width=1.5, color_U] (AM) to (Faug) ;
	
	\draw[-,>=latex, line width=1.5, color_U] (Fm) to (Caug) ;
	\draw[-,>=latex, line width=1.5, color_U] (Am) to (Caug) ;
	\draw[-,>=latex, line width=1.5, color_U] (Dbm) to (Caug) ;
	\draw[-,>=latex, line width=1.5, color_U] (CM) to (Caug) ;
	\draw[-,>=latex, line width=1.5, color_U] (AbM) to (Caug) ;
	\draw[-,>=latex, line width=1.5, color_U] (EM) to (Caug) ;

	\draw[-,>=latex, line width=1.5, color_U] (BM) to (Gaug) ;
	\draw[-,>=latex, line width=1.5, color_U] (EbM) to (Gaug) ;
	\draw[-,>=latex, line width=1.5, color_U] (GM) to (Gaug) ;
	\draw[-,>=latex, line width=1.5, color_U] (Cm) to (Gaug) ;
	\draw[-,>=latex, line width=1.5, color_U] (Em) to (Gaug) ;
	\draw[-,>=latex, line width=1.5, color_U] (Abm) to (Gaug) ;

	\draw[-,>=latex, line width=1.5, color_U] (DM) to (Daug) ;
	\draw[-,>=latex, line width=1.5, color_U] (FsM) to (Daug) ;
	\draw[-,>=latex, line width=1.5, color_U] (BbM) to (Daug) ;
	\draw[-,>=latex, line width=1.5, color_U] (Ebm) to (Daug) ;
	\draw[-,>=latex, line width=1.5, color_U] (Bm) to (Daug) ;
	\draw[-,>=latex, line width=1.5, color_U] (Gm) to (Daug) ;

	\draw[-,>=latex, line width=1.5, color_P] (Faug) to [out=140, in=220, looseness=7.0] (Faug) ;
	\draw[-,>=latex, line width=1.5, color_L, dashed] (Faug) to [out=150, in=210, looseness=4.0] (Faug) ;

	\draw[-,>=latex, line width=1.5, color_P] (Caug) to [out=50, in=130, looseness=7.0] (Caug) ;
	\draw[-,>=latex, line width=1.5, color_L, dashed] (Caug) to [out=60, in=120, looseness=6.0] (Caug) ;

	\draw[-,>=latex, line width=1.5, color_P] (Gaug) to [out=40, in=-40, looseness=7.0] (Gaug) ;
	\draw[-,>=latex, line width=1.5, color_L, dashed] (Gaug) to [out=30, in=-30, looseness=4.0] (Gaug) ;
	
	\draw[-,>=latex, line width=1.5, color_P] (Daug) to [out=230, in=310, looseness=7.0] (Daug) ;
	\draw[-,>=latex, line width=1.5, color_L, dashed] (Daug) to [out=240, in=300, looseness=6.0] (Daug) ;
		
\end{tikzpicture}
}
\end{center}
\caption{The colored Cube Dance graph for the binary relations $\mathcal{U}$ (color: black / BW: black), $\mathcal{P}$ (color: orange / BW: gray), and $\mathcal{L}$ (color: dashed green / BW: dashed gray).}
\label{fig:UPL_monoids_graphs}
\end{figure}

The `colored Cube Dance graph' (Fig. \ref{fig:UPL_monoids_graphs}) is then defined as the graph having $X$ as its set of vertices, and the set of pairs of triads related by either $\mathcal{U}$, $\mathcal{P}$, or $\mathcal{L}$ as its set of edges, each edge having a canonically attributed color in the set $\{ \mathcal{U}, \mathcal{P}, \mathcal{L} \}$. From an algebraic point of view, these binary relations generate a monoid $M_{\mathcal{U},\mathcal{P}, \mathcal{L}}$ with an action on $X$, which corresponds in categorical terms to the definition of a functor $S \colon M_{\mathcal{U},\mathcal{P}, \mathcal{L}} \to \mathbf{Rel}$, and whose structure has been investigated in \cite{Popoff2018}.

\begin{proposition}
The monoid $M_{\mathcal{U},\mathcal{P}, \mathcal{L}}$ generated by the relations $\mathcal{U}$, $\mathcal{P}$, and $\mathcal{L}$ contains 40 elements and has the following presentation.
$$\begin{aligned} M_{\mathcal{U},\mathcal{P}, \mathcal{L}} =  \langle \mathcal{U}, \mathcal{P}, \mathcal{L} \mid {} & \mathcal{P}^2=\mathcal{L}^2=e, \hspace{0.2cm} \mathcal{LPL}=\mathcal{PLP}, \hspace{0.2cm} \mathcal{U}^3=\mathcal{U}, \\
					&  \mathcal{U}\mathcal{P}=\mathcal{U}\mathcal{L}, \hspace{0.2cm} \mathcal{P}\mathcal{U}=\mathcal{L}\mathcal{U}, \hspace{0.2cm} \mathcal{U}^2\mathcal{P}\mathcal{U}^2=\mathcal{P}\mathcal{U}^2\mathcal{P}\mathcal{U}^2\mathcal{P}, \\
					&  (\mathcal{U}\mathcal{P})^2\mathcal{U}^2 = \mathcal{P} (\mathcal{U}\mathcal{P})^2\mathcal{U}^2 \mathcal{P}, \hspace{0.2cm} \mathcal{U}^2(\mathcal{P}\mathcal{U})^2 = \mathcal{P} \mathcal{U}^2(\mathcal{P}\mathcal{U})^2 \mathcal{P} \rangle \end{aligned} $$
\end{proposition}

The categorical point of view allows us to consider automorphisms of the functor $S$ (i.e. automorphisms of the monoid action), whose general definition has been given in \cite{Popoff2018,Popoff2020} and which simplifies in this case as follows.

\begin{definition}
The automorphism group $\text{Aut}(S)$ of the functor $S \colon M_{\mathcal{U},\mathcal{P}, \mathcal{L}} \to \mathbf{Rel}$ is the group of pairs $(N,\nu)$ where $N \colon M_{\mathcal{U},\mathcal{P}, \mathcal{L}} \to M_{\mathcal{U},\mathcal{P}, \mathcal{L}}$ is an automorphism, and $\nu$ is a bijection on $X$ such that we have $p \mathcal{R} q \implies \nu(p) N(\mathcal{R}) \nu(q)$ for all $\mathcal{R} \in M_{\mathcal{U},\mathcal{P}, \mathcal{L}}$ and $(p,q) \in X^2$. Composition is done term-wise.
\end{definition}

It should be noted that the normal subgroup of $\text{Aut}(S)$ of automorphisms of the form $(\text{id},\nu)$ is isomorphic to the normal subgroup of graph automorphisms leaving the color of edges invariant, some of these underlying chord progressions in pop music \cite{Popoff2020,GunnerHamiltonian}. The computation of the full automorphism group of the monoid action is therefore of interest, giving tools for musicians and composers to transform chord progressions in the colored Cube Dance. We establish its structure in the next Section.

\section{The automorphism group of the monoid action of $M_{\mathcal{U},\mathcal{P}, \mathcal{L}}$}

The automorphism group of the monoid $M_{\mathcal{UPL}}$ itself has been determined in \cite{Popoff2018}.

\begin{theorem}
The automorphism group of the $M_{\mathcal{UPL}}$ monoid is isomorphic to the group $D_6 \times \mathbb{Z}_2$.
\end{theorem}

Each automorphism $N$ of $M_{\mathcal{UPL}}$ is entirely determined by an automorphism of the subgroup isomorphic to $D_6$ generated by $\mathcal{P}$ and $\mathcal{L}$, and by the choice of the image of $\mathcal{U}$ by $N$ in the set $\{ \mathcal{U}, \mathcal{L}\mathcal{U}\mathcal{L} \}$. The main result of this paper is the structure of the automorphism group of the monoid action $S \colon M_{\mathcal{UPL}} \to \mathbf{Rel}$.

\begin{table}
	\caption{Graphical representation of the possible automorphisms $(N,\nu)$ of the functor $S \colon M_{\mathcal{U},\mathcal{P}, \mathcal{L}} \to \mathbf{Rel}$. The mapping of subsets $X_M$ is determined by the permutation of augmented chords, by the image of the generator $\mathcal{U}$ by $N$, and by the action of the group elements $g_i$ in $\mathbb{Z}_3$.}
	\label{tab:automorphisms}
	\centering
	\begin{tabular}{c|c}
		$N(\mathcal{U})=\mathcal{U}$ & $N(\mathcal{U})=\mathcal{L}\mathcal{U}\mathcal{L}$ \\
		\hline
		\autdiagrampositive{0}{3}{6}{9}{0}{3}{6}{9} &  \autdiagramnegative{0}{3}{6}{9}{0}{3}{6}{9} \\ \arrayrulecolor{lightgray}\hline \arrayrulecolor{black} \\[-0.3cm]
		\autdiagrampositive{3}{6}{9}{0}{3}{6}{9}{0} &  \autdiagramnegative{3}{6}{9}{0}{3}{6}{9}{0} \\ \arrayrulecolor{lightgray}\hline \arrayrulecolor{black} \\[-0.3cm]
		\autdiagrampositive{6}{9}{0}{3}{6}{9}{0}{3} &  \autdiagramnegative{6}{9}{0}{3}{6}{9}{0}{3} \\ \arrayrulecolor{lightgray}\hline \arrayrulecolor{black} \\[-0.3cm]
		\autdiagrampositive{9}{0}{3}{6}{9}{0}{3}{6} &  \autdiagramnegative{9}{0}{3}{6}{9}{0}{3}{6} \\ \arrayrulecolor{lightgray}\hline \arrayrulecolor{black} \\[-0.3cm]
		\autdiagramnegative{0}{9}{6}{3}{9}{6}{3}{0} &  \autdiagrampositive{0}{9}{6}{3}{9}{6}{3}{0} \\ \arrayrulecolor{lightgray}\hline \arrayrulecolor{black} \\[-0.3cm]
		\autdiagramnegative{3}{0}{9}{6}{0}{9}{6}{3} & \autdiagrampositive{3}{0}{9}{6}{0}{9}{6}{3} \\ \arrayrulecolor{lightgray}\hline \arrayrulecolor{black} \\[-0.3cm]
		\autdiagramnegative{6}{3}{0}{9}{3}{0}{9}{6} &  \autdiagrampositive{6}{3}{0}{9}{3}{0}{9}{6} \\ \arrayrulecolor{lightgray}\hline \arrayrulecolor{black} \\[-0.3cm]
		\autdiagramnegative{9}{6}{3}{0}{6}{3}{0}{9} &  \autdiagrampositive{9}{6}{3}{0}{6}{3}{0}{9} \\
	\end{tabular}
\end{table}

\begin{theorem}
The automorphism group of the functor $S \colon M_{\mathcal{UPL}} \to \mathbf{Rel}$ is a group of order 7776 isomorphic to $({\mathbb{Z}_3}^4 \rtimes D_8) \rtimes (D_6 \times \mathbb{Z}_2)$.
\end{theorem}
\begin{proof}
We sketch here the methodology for the proof, leaving the full enumeration of the cases to the reader. We denote by $\overline{C_M}$ the set $\{C_M, E_M, A\flat_M\}$, by $\overline{F_m}$ the set $\{F_m, A_m, D\flat_m\}$, and so on.

Let $N$ be an automorphism of $M_{\mathcal{U}\mathcal{P}\mathcal{L}}$. Assume for example that $N(\mathcal{U})=\mathcal{U}$. We then look for the possible bijections $\nu$ of $X$: these will obviously map the subset $\{C_{\text{aug}}, G_{\text{aug}}, D_{\text{aug}}, F_{\text{aug}}\}$ onto itself.
At this point, we can freely choose the image of $C_{\text{aug}}$ by $\nu$: assume for example that $\nu(C_{\text{aug}})=G_{\text{aug}}$.
Since $C_{\text{aug}}$ is related to the elements in $\overline{C_M} \cup \overline{F_m}$ by the relation $\mathcal{U}$, it is implied by the definition of $\nu$ and the fact that $N(\mathcal{U})=\mathcal{U}$ that $\overline{C_M} \cup \overline{F_m}$ should be bijectively mapped by to the subset $\overline{G_M} \cup \overline{C_m}$. Since the subset $\overline{C_M}$ (resp. $\overline{F_m}$) is an orbit of $C_M$ by the subgroup of $D_6$ isomorphic to $\mathbb{Z}_3$ generated by $\mathcal{L}\mathcal{P}$, we conclude by the definition of $\nu$ that either $\overline{C_M}$ is mapped to $\overline{G_M}$ and $\overline{F_m}$ to $\overline{C_m}$, or the other way around.
Assume the first case, and note that each of this mapping is entirely determined by the choice of a representative element in each subset and an element of the subgroup of $D_6$ isomorphic to $\mathbb{Z}_3$ generated by $\mathcal{L}\mathcal{P}$. We then have that $\overline{C_m}$ is mapped to $\overline{G_m}$ and $\overline{F_M}$ to $\overline{C_M}$. Since the elements of $\overline{C_m}$ are related to $G_{\text{aug}}$ by $\mathcal{U}$, it is implied by the definition of $\nu$ and the fact that $N(\mathcal{U})=\mathcal{U}$ that $G_{\text{aug}}$ should be mapped to $D_{\text{aug}}$. Similarly, we get that $F_{\text{aug}}$ should be mapped to $C_{\text{aug}}$. By continuing this enumeration, we arrive at the graphical representations of automorphisms given in Table \ref{tab:automorphisms}. The group of permutations of the set $\{C_{\text{aug}}, G_{\text{aug}}, D_{\text{aug}}, F_{\text{aug}}\}$ is isomorphic to $D_8$, and this automatically determines the permutation of the set of subsets $\{ \overline{C_M}, \overline{G_M}, \overline{D_M}, \overline{F_M} \}$. For a given permutation of this set, each element $g_i$ is a group element in $\mathbb{Z}_3$ determining how the subsets are mapped, assuming a set of representative elements has been fixed beforehand. It can then readily be seen that $\text{Aut}(S)$ is isomorphic to $({\mathbb{Z}_3}^4 \rtimes D_8) \rtimes (D_6 \times \mathbb{Z}_2)$, a group of order 7776.
\end{proof}

\section{An interactive interface for composing (with) automorphisms}

\begin{figure}[b!]
	\begin{center}
	\includegraphics[scale=0.24]{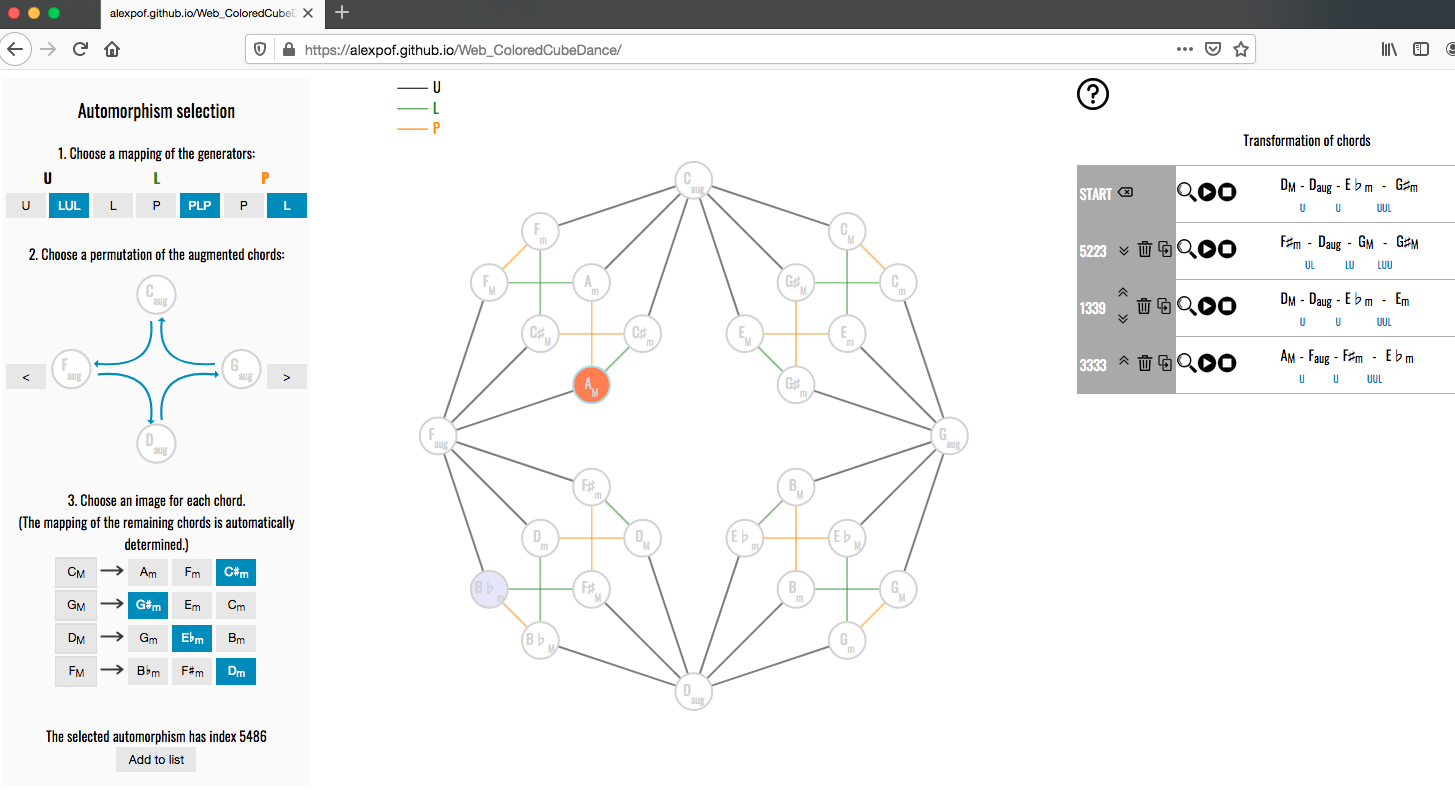}
	\end{center}
	\caption{Screenshot of the web interface for manipulating the automorphisms of the colored Cube Dance.}
	\label{fig:interface_screenshot}
\end{figure}

Contrary to the neo-Riemannian $PLR$ group, the size and complexity of $\text{Aut}(S)$ makes it hard to use with pen and paper, especially for non-mathematicians wishing to explore its potential for transforming chord progressions. In this light, we have developed a web interface for the concrete manipulation of the elements of $\text{Aut}(S)$, and for their use as chord transformations. A screenshot of this interface is shown in Figure \ref{fig:interface_screenshot}. It uses common HTML and Javascript frameworks \cite{D3JS}, thus making it runnable on virtually any web browser without the need for complicated software installations. Such frameworks have already been used for other mathematics/music applications, notably to explore the Tonnetz \cite{TonnetzGuichaoua}. The web interface can directly be used from the corresponding GitHub repository \cite{WebColorCubeDance} and the associated source code is freely available.

As shown on Figure \ref{fig:interface_screenshot}, the left part of the interface corresponds to the interactive choice of an automorphism of the colored Cube Dance. The middle part is an interactive colored Cube Dance: alt-clicking on chords adds them to the current chord progression, which is shown on the right part of the interface, along with its successive transformation by the selected automorphisms. Each chord progression can be played back with the corresponding buttons. In the future, the interface will also feature MIDI capabilities, so that chord progressions could be recorded, transformed, and replayed at will (see \cite{TonnetzGuichaoua} for a current implementation of such capabilities). Following the pattern of Table \ref{tab:automorphisms}, the user first selects a mapping of the generators of $M_{\mathcal{UPL}}$, then a permutation of the augmented chords, and finally a mapping of the major/minor chords through the mapping of a given representative in each quadrant. Once an automorphism has been uniquely determined, the user can hover over chords in the middle representation of the colored Cube Dance to see how they are mapped by the selected element of $\text{Aut}(S)$. The `add to list' button appends the selected automorphism to the list on the right, in which the current chord progression is successively transformed through automorphism composition.

The combination of SVG graphics possibilities in HTML with Javascript allows one to quickly develop user-friendly interfaces for math/music concepts, thus showing the potential of these technologies for outreach activities. It is our hope that the colored Cube Dance web interface will prove useful for students, musicians, and composers to creatively explore chord transformations \text{via} automorphisms.

\bibliography{coloredcubedance_automorphisms}
\bibliographystyle{alpha}
\addcontentsline{toc}{section}{References}

\end{document}